\documentclass[12pt,
amscd]{amsart}
\usepackage{latexsym,amssymb}  
\theoremstyle{plain}
  \newtheorem{thm}{Theorem}[section]
  \newtheorem{prop}[thm]{Proposition}
  
  \newtheorem{cor}[thm]{Corollary}
\theoremstyle{definition}
  \newtheorem{dfn}[thm]{Definition}
  
\theoremstyle{remark}
  \newtheorem{rem}[thm]{Remark}
\numberwithin{equation}{section}
\def\ba{{\bf a}}
\def\bb{{\bf b}}
\def\b1{{\mathbf 1}}
\def\bx{{\mathbf x}}
\def\rH{\tilde{H}}
\def\rh{\tilde{h}}

\def\M{M^\bullet}
\def\N{N^\bullet}
\def\P{P^\bullet}
\def\PP{{\mathbb P}}
\def\I{I^\bullet}

\def\cS{{\mathcal S}}
\def\cE{{\mathcal E}}
\def\bA{{\mathbf {A}}}

\def\p{{\mathfrak p}}
\def\NN{{\mathbb N}}
\def\Sq{\operatorname{Sq}_S}
\def\SqE{\operatorname{Sq}_E}
\def\ZZ{{\mathbb Z}}

\def\bL{{\mathbf L}}
\def\cL{{\mathcal L}}
\def\bR{{\mathbf R}}
\def\Hom{\operatorname{Hom}}
\def\Spec{\operatorname{Spec}}
\def\Ext{\operatorname{Ext}}
\def\Coh{\operatorname{Coh}}
\def\Tor{\operatorname{Tor}}
\def\Top{\operatorname{Top}}
\def\Sym{\operatorname{Sym}}
\def\Ann{\operatorname{Ann}}

\def\ann{\operatorname{Ann}}
\def\cone{\operatorname{cone}}
\def\supp{\operatorname{supp}}
\def\reg{\operatorname{reg}}
\def\<{{\langle}}
\def\>{{\rangle}}

\def\pd{\operatorname{proj.dim}}
\def\D{\omega^\bullet}
\def\MDS{\operatorname{Gr} S}
\def\mdE{\operatorname{gr} E}
\def\mod{\operatorname{gr}}
\def\mdS{\operatorname{gr} S}

\def\DS{{\bf D}_S}
\def\DE{{\bf D}_E}
\def\mmEn{\operatorname{*gr} E}
\def\mmSn{\operatorname{*gr} S}
\begin{document}
\title[BGG correspondence and R\"omer's theorem on an exterior algebra]
{BGG correspondence and R\"omer's theorem\\ on an exterior algebra}
\dedicatory{To the memory of Professor Tetsushi Ogoma}
\author{Kohji Yanagawa}
\address{Department of Mathematics, 
Graduate School of Science, Osaka University, Toyonaka, Osaka 
560-0043, Japan}
\email{yanagawa@math.sci.osaka-u.ac.jp}

\begin{abstract}
Let $E = K\< y_1, \ldots, y_n \> $ be the exterior algebra. 
The {\it (cohomological) distinguished pairs} of a graded $E$-module $N$
describe the growth of a minimal graded injective resolution of $N$. 
R\"omer gave a duality theorem between the distinguished pairs of $N$ 
and those of its dual $N^*$. In this paper, we show that 
under Bernstein-Gel'fand-Gel'fand correspondence, his theorem is 
translated into a natural corollary of 
local duality for (complexes of) graded 
$S=K[x_1, \ldots, x_n]$-modules. Using this idea, 
we also give a $\ZZ^n$-graded version of R\"omer's theorem. 
\end{abstract}

\maketitle
\section*{Introduction}
In this section, to introduce a background of the present paper, 
we summarize results of Aramova-Herzog~\cite{AH} and R\"omer~\cite{R0}. 

Let $S = K[x_1, \ldots, x_n]$ be the polynomial ring over a field $K$, and 
$M$ a finitely generated graded $S$-module. The $ij$th Betti number 
$\beta_{i,j}(M) = \dim_K \Tor^S_i(K,M)_j$ of $M$ 
is an important invariant. Following Bayer-Charalambous-Popescu~\cite{BCP}, 
we say a Betti number $\beta_{k,m}(M) \ne 0$ is {\it extremal}, 
if $\beta_{i,j}(M) = 0$ for all $(i,j) \ne (k,m)$ with $i \geq k$ and 
$j-i \geq m-k$. This notion has two remarkable properties. 
First, a homogeneous ideal $I \subset S$ has the same extremal Betti numbers 
as its generic initial ideal $\operatorname{Gin}(I)$ with respect to the 
reverse lexicographic term order on $S$. 
Another important property is the following. 

\medskip

\noindent{\bf Theorem~A} 
(Bayer-Charalambous-Popescu, \cite[Theorem~2.8]{BCP}) \ 
{\it Let $\Delta \subset 2^{\{1, \ldots, n\}}$ be a simplicial complex, and 
$K[\Delta] = S/I_\Delta$ the Stanley-Reisner ring. 
And let $\Delta^\vee$ be the Alexander dual complex of $\Delta$. 
Then $\beta_{i, i+j}(K[\Delta])$ is  extremal if and only if 
so is $\beta_{j, i+j}(I_{\Delta^\vee})$. 
Moreover, if this is the case,  
then $\beta_{i, i+j}(K[\Delta]) = \beta_{j, i+j}(I_{\Delta^\vee})$. }

\medskip

We have $\beta_{i, n}(K[\Delta]) = \dim_K \rH_{n-i-1}(|\Delta|;K) 
=: \rh_{n-i-1}(|\Delta|)$ 
by Hochster's formula. If $\beta_{i, n}(K[\Delta]) \ne 0$ then 
it is always an extremal Betti number. The equality 
$\rh_{n-i-1}(|\Delta|) = \beta_{i, n}(K[\Delta]) =  
\beta_{n-i,n}(I_{\Delta^\vee}) = \beta_{n-i+1,n}(K[\Delta^\vee]) = 
\rh_{i-2}(|\Delta^\vee|)$ induced by Theorem~A corresponds to 
the usual Alexander duality. More generally, 
Theorem~A gives an Alexander duality for (some) {\it iterated Betti numbers} 
(c.f. \cite{Kalai, BTN}). 

Let $E = K\< y_1, \ldots, y_n \>$ be the exterior algebra. 
To understand Theorem~A, Aramova-Herzog~\cite{AH} introduced 
{\it distinguished pairs} for a graded $E$-module $N$. 
See Definition~\ref{dist pair for E} below. 
({\it We use a different convention to describe 
these pairs. See Remark~\ref{convention}.})  
The distinguished pairs of $N$ roughly describe the growth of 
the minimal graded (infinite) injective resolution of $N$.
Let $K\{\Delta\} = E/J_\Delta$ be the exterior face ring of $\Delta$. 
Then \cite[Corollary~9.6]{AH} states that $(d,i)$ is a distinguished 
pair for $K\{\Delta\}^* := \Hom_E(K\{\Delta\}, E)$ 
if and only if $\beta_{d+i-n,d}(K[\Delta])$ is extremal. 
And the values of extremal Betti numbers 
can be described by the Cartan cohomologies of $K\{\Delta\}^*$.

R\"omer proved that  $(d,i)$ is distinguished for $N$ if and only if 
so is $(d, 2n-d-i)$ for $N^*$. He also gave a duality between 
certain components of Cartan cohomologies of $N$ and those of $N^*$. 
Since $k\{\Delta\}^* = J_{\Delta^\vee}$, his result implies Theorem~A.  

Bernstein-Gel'fand-Gel'fand correspondence (BGG correspondence, for short) 
is a well known theorem which states that 
the derived category $D^b(\mdS)$ of finitely generated graded 
$S$-modules is equivalent to the similar category $D^b(\mdE)$ for $E$. 
In this paper, we give a new proof of the result of R\"omer 
using BGG correspondence. 
More precisely, under this correspondence, R\"omer's theorem 
is translated into a statement on $D^b(\mdS)$ which is a natural 
consequence of the local duality (Serre duality). 
A key point is that the duality 
functor $\Hom_E(-,E)$ on $D^b(\mdE)$ corresponds to the duality functor 
$R\Hom_S(-,\D)$ on $D^b(\mdS)$, where $\D$ is a dualizing complex of $S$. 
So our proof is (philosophically) 
simple and sheds new light on the result. 

The original paper \cite{BCP} states Theorem~A in the 
$\ZZ^n$-graded context, while the arguments in \cite{AH, R0} 
are hard to work in this context. 
But, since BGG correspondence also holds for $\ZZ^n$-graded modules, 
our method is powerful in this context too. 
See \S2.  This part of the present paper is a continuation 
of the author's previous paper \cite{Y1}. 

\section{$\ZZ$-graded case}
Let $W$ be an $n$-dimensional vector space over a field $K$, and 
$S = \bigoplus_{i \geq 0} \Sym_i W$ the polynomial ring.  
We regard $S$ as a graded ring with $S_i = \Sym_i W$.
Let $\MDS$ be the category of graded $S$-modules and their degree preserving 
$S$-homomorphisms, and $\mdS$ the full subcategory of $\MDS$ consisting of 
finitely generated modules. Then there is an
equivalence $D^b(\mdS) \cong D^b_{\mdS}(\MDS)$. 
(For derived categories, consult \cite{GM}.)
So we will freely identify these categories. 
For $M = \bigoplus_{i \in \ZZ} M_i \in \MDS$ 
and an integer $j$, $M(j)$ denotes the shifted module with $M(j)_i = M_{i+j}$. 
For $\M \in D^b(\MDS)$, $\M[j]$ denotes the $j$th translation of $\M$, 
that is, $\M[j]$ is the complex with $\M[j]^i = M^{i+j}$. 
So, if $M \in \MDS$, $M[j]$ is 
the cochain complex $\cdots \to 0 \to M \to 0 \to \cdots$, 
where $M$ sits in the  $(-j)$th position. 
If $M \in \mdS$ and $N \in \MDS$, then $\Hom_S(M,N)$ 
has the structure of a graded $S$-module with  
$\Hom_S(M,N)_i = \Hom_{\MDS}(M,N(i))$. 

Let $\D \in D^b(\mdS)$ be a minimal graded injective resolution of $S(-n)[n]$. 
That is, $\D$ is a graded normalized dualizing complex of $S$. 
Then $\DS(-) := \Hom^\bullet_S(-, \D)$ gives a 
duality functor from $D^b(\mdS)$ to itself. The $i$th cohomology of 
$\DS(\M)$ is $\Ext^i_S(\M,\D)$. For $\M \in D^b(\mdS)$ and $i \in \ZZ$, set 
$d_i (\M):= \dim_S H^i(\M)$. Here the Krull dimension of the 0 module 
is $-\infty$. 

\begin{dfn}\label{dist pair for S}
We say $(d, i) \in \NN \times \ZZ$ is a 
{\it distinguished pair} for a complex $\M \in D^b(\mdS)$, 
if $d = d_i(\M)$ and 
$d_j(\M) < d + i- j$ for all $j$ with $j < i$. 
\end{dfn}

Let $\M \in D^b(\mdS)$ and $d = d_i(\M) \geq 0$. 
If $d = \max\{ \, d_j(\M) \mid j \in \ZZ \, \}$, 
then $(d,i)$ is distinguished for $\M$. 
On the other hand, if $i = \min\{ \, j \mid H^j(\M) \ne 0 \, \}$, 
then $(d,i)$ is also distinguished. 
Thus $\M$ has several distinguished pairs in general. 

In this paper, $\deg_S(M)$ denotes the multiplicity of a module $M \in \mdS$ 
(i.e., $e(M)$ of \cite[Definition~4.1.5]{BH}).

\begin{thm}\label{main}
For $\M \in D^b(\mdS)$, we have the following. 

(1)  A pair $(d, i)$ is distinguished for $\M$ 
if and only if $(d, -d-i)$ is distinguished for $\DS(\M)$. 

(2) If $(d, i)$ is a distinguished pair for $\M$, 
then $$\deg_S H^i(\M) = \deg_S \Ext_S^{-d-i}(\M, \D).$$  
\end{thm}

\begin{proof}
(1) Since the statement is ``symmetric", it suffices to 
prove the direction $\Rightarrow$. 

From the double complex $\Hom^\bullet_S (\M, \D)$, we have 
a spectral sequence $E_2^{p,q} = \Ext^p_S(H^{-q}(\M), \D) 
\Rightarrow \Ext^{p+q}_S(\M, \D)$.  
For simplicity, set $e_r^{p,q} := \dim_S E_r^{p,q}$. 
Since $\Ext^i_S(M, \D) \cong  \Ext^{n+i}_S(M,S(-n))$ for $M \in \mdS$, 
the following inequality follows from argument analogous to  
\cite[\S 8.1, Theorem~8.1.1]{BH}.
\begin{equation}
e_2^{p,q} = \dim_S \Ext^p_S(H^{-q}(\M), \D) = 
\begin{cases}\label{basic}
-p &  \text{if $p = -d_{-q}(\M)$,}\\
\leq -p & \text{if $-d_{-q}(\M) < p \leq 0$,}\\
- \infty & \text{otherwise.}
\end{cases}
\end{equation}

(I) \ By \eqref{basic}, we have $e_2^{-d, -i} = d$.  
On the other hand, we have $e_2^{p,q} < d$ for all 
$(p,q) \ne (-d, -i)$ with $p+q = -d-i$. In fact, 
the assertion follows from \eqref{basic} if $p > -d$. 
So we may assume that $p < -d$ and $q = -d-i-p > -i$. 
Since $(d, i)$ is distinguished, we have 
$d_{-q} (\M) < d + i + q = -p$. Thus $E_2^{p,q} = 0$ in  this case. 
Anyway, we have $e_\infty^{p,q} < d$ for all 
$(p,q) \ne (-d, -i)$ with $p+q = -d-i$. 

\smallskip

(II) \ Since $d_{i-j+1}(\M) < d + j-1 < d+j$ for all $j \geq 2$, 
we have that $E_2^{-d-j, -i+j-1} = 0$. So we have 
$E_r^{-d-j, -i+j-1} = 0$ for all $r \geq 2$. 
Next we will show that $d = e_2^{-d, -i} =  e_3^{-d, -i} = \cdots 
=  e_r^{-d, -i}$ by induction on $r$. 
Recall that $E_{r+1}^{-d, -i}$ is the cohomology of 
$$E_r^{-d-r, -i+r-1} \to E_r^{-d, -i} \to E_r^{-d+r, -i-r+1}.$$
But we have seen that $E_r^{-d-r, -i+r-1} = 0$. 
Moreover,  $e_r^{-d+r, -i-r+1} \leq e_2^{-d+r, -i-r+1} 
\leq  d-r < d$ by \eqref{basic}, 
and $e_r^{-d,-i}= d$ by the induction hypothesis. 
Thus $e_{r+1}^{-d,-i}= d$. Hence $e_{\infty}^{-d,-i}= d$. 
From this fact and (I), 
we have that $\dim_S \Ext^{-d-i}_S(\M, \D) = d$. 

\smallskip

(III) \ Finally, we will show that $\dim_S \Ext^{-d-i-j}_S(\M, \D) < d+j$ 
for all $j > 0$. 
To see this, it suffices to show that $e_2^{p,q} < d+j$ for all 
$j > 0$ and all $(p, q)$ with $p+q = -d-i-j$. If $p > -d-j$, the assertion 
is clear. If $p = -d-j$, then $q = -i$ and $d_{-q}(\M) = d < -p$. So 
$E_2^{p,q} = 0$ in this case. Hence we may assume that $p < -d -j$ and 
$-q =d+i+j+p < i$. Since $(d, i)$ is distinguished, $d_{-q}(\M) <  
d+(i+q) = -j-p < -p$. So we have $E_2^{p,q} = 0$ in this case too. 

\smallskip

(2) Since $\deg_S E_r^{-d, -i} = \deg_S E_{r+1}^{-d,-i}$ for all $r \geq 2$ 
by the argument in (II) of the proof of (1), we have 
$\deg_S E_2^{-d, -i} = \deg_S E_\infty^{-d,-i}$. 
So we have 
$$\deg_S \Ext_S^{-d-i}(\M, \D) = \deg_S E_\infty^{-d, -i} = 
 \deg_S E_2^{-d, -i} = \deg_S \Ext_S^{-d}(H^i(\M), \D),$$
where the first equality follows from (I) and (II). 

For a module $M \in \mdS$ of dimension $d$, we have 
$\deg_S M = \deg_S \Ext_S^{-d}(M, \D)$. In fact, for a prime ideal 
$\p \subset S$ with $\dim S/\p = d$, let $\bx$ be a maximal $S_\p$-sequence 
contained in $\ann_{S_\p}(M_\p)$, and $R := S_\p/\bx S_\p$ an artinian 
Gorenstein local ring. Then we have  
$$\Ext_S^{-d}(M, \D) \otimes_S S_{\p} \cong \Ext_{S_\p}^{n-d}(M_\p, S_\p) 
\cong \Hom_R(M_\p, R).$$
Therefore, $$l_{S_p}(\Ext_S^{-d}(M, \D) \otimes_S S_\p) 
= l_R(\Hom_R(M_\p, R)) = l_R(M_\p) = l_{S_p}(M_\p).$$ 

Since $\dim_S (H^i(\M)) =d$, we have 
$\deg_S \Ext_S^{-d}(H^i(M), \D) = \deg_S H^i(\M)$. 
\end{proof}

\begin{rem}
Theorem~\ref{main} (1) also holds for a noetherian local ring $A$
admitting a dualizing complex. The part (2) also holds for $A$, 
if we replace $\deg_S (-)$ by $l_{A_\p} (- \otimes_A {A_\p})$ 
for a prime ideal $\p \subset A$ with $\dim A/\p = d$. 
\end{rem}

Let $V$ be the dual vector space of $W$, and $E = \bigwedge V$ the 
exterior algebra. 
We regard $E$ as a negatively graded ring with 
$E_{-i} = \bigwedge^i V$ (this is the opposite convention from 
\cite{AH, R0}). Let $\mdE$ be the category of finitely generated 
graded $E$-modules and their degree preserving $E$-homomorphisms. 
Here ``$E$-module" means a left and right module $N$ with 
$ea= (-1)^{(\deg e)(\deg a)}ae$ for all homogeneous elements $e \in E$ and 
$a \in N$. 

Let $\{x_1, \ldots, x_n \}$ be a basis of $W$, and $\{y_1, \ldots, y_n\}$ 
its dual basis of $V$. For a complex $\N$ in $\mdE$, 
set $\bL(\N) = \bigoplus_{i \in \ZZ} S \otimes_K N^i$ 
and $\bL(\N)^m = \bigoplus_{i-j=m} S \otimes_K N_j^i$.  
The differential defined by 
$$\bL(\N)^m \supset S \otimes_K N_j^i \ni 1 \otimes z \mapsto 
\sum_{1 \leq l \leq n} x_l  \otimes y_l z + (-1)^m (1 \otimes \delta^i(z)) 
\in \bL(\N)^{m+1}$$
makes $\bL(\N)$ a cochain complex of free $S$-modules. 
Here $\delta^i$ is the $i$th differential map of $\N$. 
Moreover, $\bL$ gives a functor from $D^b(\mdE)$ to $D^b(\mdS)$. 

For $M \in \mdS$ and $i \in \ZZ$, we can define a graded $E$-module 
structure on $\Hom_K(E, M_i)$ by $(a f)(e) = f(ea)$. 
Then $\Hom_K(E, M_i) \cong E(-n) \otimes_K M_i$.  
Set $\bR(M) = \Hom_K(E, M)$ and $\bR^i(M) = \Hom_K(E, M_i)$. 
The differential defined by 
$$\bR^i(M) = \Hom_K(E,M_i) \ni f \mapsto 
[e \mapsto \sum_{1 \leq j \leq n} 
x_j f(y_je)] \in \Hom_K(E, M_{i+1}) = \bR^{i+1}(M)$$
makes $\bR(M)$ a cochain complex of free $E$-modules.
We can also construct $\bR(\M)$ from a complex $\M$ in natural way. 
Then $\bR$ gives a functor from 
$D^b(\mdS)$ to $D^b(\mdE)$. See \cite{EFS} for details.  
The following is a crucial result. 

\begin{thm}[BGG correspondence, c.f.{\cite{EFS}}]
The functors $\bL$ and $\bR$ give a category equivalence 
$D^b(\mdS) \cong D^b(\mdE)$. 
\end{thm}

For $N \in \mdE$, then $N^* := \Hom_E(N,E) \cong \Hom_K(N,K)(n)$ 
is a graded $E$-module again. 
$(-)^*$ gives an exact duality functor on $\mdE$, and it can be extended 
to the duality functor $\DE$ on $D^b(\mdE)$. 

\begin{prop}\label{dualities}
For $\N \in D^b(\mdE)$, we have $$\DS \circ \bL (\N) \cong 
\bL \circ \DE (\N) (-2n)[2n].$$ 
\end{prop}

\begin{proof}
Since $\bL(\N)$ consists of free $S$-modules, we have 
$$\DS \circ \bL(\N) \cong \Hom_S^\bullet(\bL(\N), S(-n)[n]).$$ 
It is easy to see that 
$$\Hom_S^m(\bL(\N), S(-n)[n]) \cong \bigoplus_{j-i=m+n} 
S(-n) \otimes_K (N_j^i)^\vee,$$
where $(-)^\vee$ means the graded $K$-dual. 
On the other hand, 
\begin{eqnarray*}
\bL \circ \DE(\N)^m 
= \bigoplus_{i-j=m} S \otimes_K \DE(\N)_j^i 
&=& \bigoplus_{i-j=m} S(n) \otimes_K (N_{-n-j}^{-i})^\vee \\
&=& \bigoplus_{j-i=m-n} S(n) \otimes_K (N_j^i)^\vee.
\end{eqnarray*}
So we can easily construct a quasi-isomorphism 
$\DS \circ \bL (\N) \to 
\bL \circ \DE (\N) (-2n)[2n].$ 
\end{proof}

For $\N \in D^b(\mdE)$, we have 
$H^i(\bL(\N))_j \cong \Ext_E^{j+i}(K, \N)_j$ by \cite[Theorem~3.7]{EFS}.  
So the Laurent series
$P_i(t)=\sum_{j \in \ZZ} (\dim_K \Ext_E^{j+i}(K, \N)_j) \cdot  t^j$ 
is the Hilbert series of the finitely generated 
graded $S$-module $H^i(\bL(\N))$. If 
$H^i(\bL(\N)) \ne 0$, there exists a Laurent polynomial 
$Q_i(t)\in \ZZ[t, t^{-1}]$ such that 
$$P_i (t) =\frac{Q_i(t)}{(1-t)^d},$$ 
where $d = d_i(\bL(\N)) = \dim_S H^i(\bL(\N))$. 
Set $e_i(\N) := Q_i(1) = \deg_S H^i(\bL(\N))$. 
So $d_i(\bL(\N))$ and $e_i(\N)$ measure the growth of the 
``$(-i)$-linear strand" of a minimal injective resolution of $\N$.  

\cite{AH, R0} treated  $d_i(\bL(N))$ and $e_i(N)$ for a module $N \in \mdE$ 
more or less indirectly. 
But their approach is very different from ours. 
They use {\it Cartan (co)homology} of $N$. 
See \cite{AH, R0} for the definition of this (co)homology. 
Let ${\bf v} = v_1, \ldots, v_n$ be a basis of $V$ which is {\it generic} 
with respect to $N$ in the sense of \cite[Definition~4.7]{AH}. 
As \cite{AH, R0}, we set $H_i(k) := H_i(v_1, \ldots, v_k;N)$ and 
$H^i(k) := H^i(v_1, \ldots, v_k;N)$ to be Cartan (co)homologies. 
Note that $H_i(n) = \Tor_i^E(K,N)$, $H^i(n) = \Ext^i_E(K,N)$ and 
$H^i(v_1, \ldots, v_k;N^*) \cong 
H_i(v_1, \ldots, v_k;N)^*$.  
It follows from the argument in \S6 of \cite{AH} that  
the function $j \mapsto \dim_K H^{i+j}(k)_j$ is a polynomial function 
for $j \gg 0$. Moreover, $d_i(\bL(N)) \leq 0$ if and only if 
$H^{i+j}(n)_{j} =0$ for $j \gg 0$ if and only if 
$H^{i+j}(k)_{j} =0$ for all $k \leq n$ and $j \gg 0$. 
\cite[Proposition~9.4]{AH} can be restated as follows: 
If $d_i(\bL(N)) > 0$, we have 
\begin{equation}\label{AB formula}
d_i(\bL(N)) = n+1 - \min\{ \, k \mid \text{$H^{i+j}(k)_j \ne 0$ 
for all $j \gg 0 $} \, \}.
\end{equation} 

A ({\it cohomological}) {\it distinguished pair} for a module 
$N \in \mdE$ was introduced in \cite[Definition~3.4]{R0} 
(see also \cite{AH}). Here we generalize this notion to 
a complex.  

\begin{dfn}\label{dist pair for E}
Let $\N \in D^b(\mdE)$. We say $(d,i) \in \NN \times \ZZ$ 
is a {\it distinguished pair} for $\N$ 
if and only if it is distinguished for $\bL(\N)$ (in the sense of 
Definition~\ref{dist pair for S}). 
\end{dfn}

\begin{rem}\label{convention}
By \eqref{AB formula}, we see that $(d,i)$ is a distinguished pair for 
a module $N \in \mdE$ in the above sense if and only if $(n+1-d,i)$ is a 
``cohomological distinguished pair" for $N$ in the sense of \cite{R0}. 
(Recall that $E$ is a positively graded ring in \cite{AH,R0}.) 
\cite{AH} also use the term ``distinguished pair". But this 
is ``homological distinguished pair" of \cite{R0}, and 
$(d,i)$ is a distinguished pair for $N$ in our sense 
if and only if $(n+1-d, n-i)$ is a distinguished pair for $N^*$ 
in the sense of \cite{AH}.
\end{rem}

\begin{cor}[c.f. {\cite[Theorem~3.8]{R0}}]\label{dist dual}
Let $\N \in D^b(\mdE)$. A pair $(d,i)$ is distinguished for $\N$ 
if and only if $(d, 2n-d-i)$ is distinguished for $\DE(\N)$. 
If this is the case, we have $e_i(\N) = e_{2n-d-i}(\DE(\N))$. 
\end{cor}

\begin{proof}
For the first statement, it suffices to prove the direction $\Rightarrow$. 
By Theorem~\ref{main}, $(d, -d-i)$ is a distinguished pair for 
$\DS \circ \bL (\N) \cong \bL \circ \DE (\N) (-2n)[2n]$. 
For a complex $\M \in D^b(\mdS)$, we have 
$H^j(\M(-2n)[2n]) = H^{2n+j}(\M)(-2n)$ and $d_j(\M(-2n)[2n]) = 
d_{2n+j}(\M)$. Thus $(d, 2n-d-i)$ is distinguished for 
$\bL \circ \DE (\N)$. The last equality follows from Theorem~\ref{main} (2). 
\end{proof}

For a module $N \in \mdE$, $d_i(\bL(N))$ can be 0 quite often. 
But we have the following. 

\begin{prop}
Assume that a module $N \in \mdE$ does not have a free summand. 
If $(d,i)$ is a distinguished pair for $N$, 
then we have $d > 0$. 
\end{prop}

\begin{proof}
Let $0 \to N \to I^0 \to I^1 \to \cdots$ (resp. $\cdots \to I^{-1} \to 
I^0 \to N \to 0$) be a minimal injective (resp. projective) 
resolution of $N$. For $j \geq 0$, set 
 $\Omega_j(N) := (\ker (I^j \to I^{j+1}))[-j]$. 
Obviously, $0 \to \Omega_j(N) \to I^j \to I^{j+1} \to \cdots$ 
is a  minimal injective resolution. On the other hand, since 
$N$ does not have a free summand, $\cdots \to I^{-1} \to 
I^0 \to \cdots \to I^{j-1} \to  \Omega_j(N) \to 0$ 
is a minimal projective resolution. If $d_i(\bL(N)) > 0$, then 
$d_i(\bL(\Omega_j(N))) = d_i(\bL(N))$ for all $j \geq 0$. But, 
if $d_i(\bL(N)) = 0$, then $d_i(\bL(\Omega_j(N))) = -\infty$ for $j \gg 0$. 
On the other hand, since a minimal injective resolution of 
$N^*$ is the dual of a minimal projective resolution of $N$, 
we have $d_i(\bL(N^*)) = d_i(\bL(\Omega_j(N)^*))$ for all $i$ and 
all $j \geq 0$. So $N^*$ and $\Omega_j(N)^*$ have the same 
distinguished pairs. 
For a contradiction, we assume that $(0,i)$ is a distinguished pair for $N$. 
Then $(0, 2n-i)$ is a distinguished pair for $N^*$ and 
$\Omega_j(N)^*$. So $(0,i)$ is a distinguished pair for 
$\Omega_j(N)$ for all $j \geq 0$. This contradicts the above observation. 
\end{proof}

We say a distinguished pair $(d,i)$ is {\it positive}, if $d > 0$. 
Since \cite{AH,R0} study a distinguished pair for a module, 
they only treat a positive one. 

\begin{rem}\label{complex -> module}
When $\N$ is a module, Corollary~\ref{dist dual} was proved in 
\cite[Theorem~3.8]{R0}. On the other hand, for positive distinguished pairs, 
we can prove the corollary from \cite[Theorem~3.8]{R0} directly:  
Let $\I$ be an injective resolution of $\N$ and $\P$ a projective 
resolution of $\I$. 
From the quasi-isomorphism $f: \P \to \I$, we have the exact complex 
$(T^\bullet, \partial^\bullet) := \cone(f)$. Then $N := \ker \partial_0$ 
(resp. $N^*$) has the same positive distinguished pairs as $\N$ 
(resp. $\DE(\N)$). 

A variant of BGG correspondence gives an equivalence 
$\underline{\mod} E \cong D^b(\Coh(\PP^{n-1}))$ of triangulated categories, 
where $\underline{\mod} \, E$ is the stable category, 
and $\Coh(\PP^{n-1})$ is the category of coherent sheaves on 
$\PP^{n-1} = \operatorname{Proj} S$. More precisely, the composition of 
the functor $\bL: \mdE \to D^b(\mdS)$ and the natural functor 
$D^b(\mdS) \to D^b(\Coh(\PP^{n-1}))$ induces this equivalence.
Note that the functor 
$\mdS \ni M \to \tilde{M} \in \Coh(\PP^{n-1})$ ignores modules 
of finite length. Hence if $d_i(\M) = 0$ then $H^i(\tilde{M}^\bullet) =0$. 
In this sense, the duality in \cite{R0} corresponds to a duality on 
$D^b(\Coh(\PP^{n-1}))$. 
\end{rem}

In the rest of this section, we assume that $K$ is algebraically closed.
Let $N \in \mdE$. Following \cite{AAH}, 
we say $v \in E_{-1} = V$ is $N$-{\it regular} if $\Ann_N(v) = vN$. 
It is easy to see that $v$ is $N$-regular 
if and only if it is $N^*$-regular. We say 
$V_E(N) = \{ v \in V \mid \text{$v$ is {\it not} $N$-regular}\}$ 
is the {\it rank variety} of $N$ (see \cite{AAH}). \cite[Theorem~3.1]{AAH} 
states that $V_E(N)$ is an algebraic subset of $V = \Spec S$, and 
$\dim V_E(V) = \max\{ \,  d_i(\bL(N)) \mid i \in \ZZ  \, \}$. 
By the above remark, $V_E(N) = V_E(N^*)$. 
We can refine this observation using the grading of $N$. 

Recall that $S$ can be seen as the Yoneda algebra $\Ext^*_E(K,K)$, and 
$\Ext^*_E(K,N)$ has the $S$-module structure. 
By the same argument as \cite[Theorem~3.9]{AAH} (see also the proof 
 of \cite[Corollary~3.2 (b)]{EPY}),  we have that 
$$V_E(N) = \{ \, v \in V \mid \text{$\xi(v) = 0$ for all 
$\xi \in \Ann_S( \, \Ext_E^*(K,N) \, )$} \, \}.$$
But $[\Ext_E^{*+i}(K,N)]_* := \bigoplus_{j \in \ZZ} 
\Ext_E^{j+i}(K,N)_j$ is an $S$-module which is isomorphic to 
$H^i(\bL(N))$ (see the proof of \cite[Proposition~2.3]{EFS}), 
and we have
$\Ext_E^*(K,N) \cong \bigoplus_{j \in \ZZ} [\Ext_E^{*+i}(K,N)]_*$. Set 
$$V_E^i(N) = \{ \, v \in V \mid \text{$\xi(v) = 0$ for all 
$\xi \in \Ann_S( \, [\Ext_E^{*+i}(K,N)]_* \, )$} \, \}.$$
We have $V_E(N) = \bigcup_i V_E^i(N)$ and 
$d_i(\bL(N)) = \dim V_E^i(N)$. For an algebraic set $X \subset \Spec S$ 
of dimension $d$, set $\Top (X)$ to be the union of the all irreducible 
components of $X$ of dimensions $d$.  

\begin{prop}
If $(d, i)$ is a distinguished pair for $N \in \mdE$, then we have 
$\Top(V_E^i(N)) = \Top(V_E^{2n-d-i}(N^*))$. 
\end{prop}

\begin{proof}
By the proof of Theorem~\ref{main}, $\Ann_S(\, H^i(\bL(N)) \,)$ 
has the same top dimensional components
as $\Ann_S(\, H^{-d-i}(\DS \circ \bL(N))\,)$. 
\end{proof}

In the above situation, we have $V_E^i(N) \ne V_E^{2n-d-i}(N^*)$ in general.

\section{Squarefree case}
In this section, we regard $S=K[x_1, \ldots, x_n]$ as an $\NN^n$-graded ring 
with $\deg x_i = (0, \ldots, 0,1,0, \ldots, 0)$ where 1 is in the $i$th 
position. Similarly, $E = K\<y_1, \ldots, y_n\>$ is a $-\NN^n$-graded ring 
with $\deg y_i = - \deg x_i$. Let $\mmSn$ (resp. $\mmEn$) be the category 
of finitely generated $\ZZ^n$-graded $S$-modules (resp. $E$-modules). 
The functors $\bL$ and $\bR$ defining the BGG correspondence $D^b(\mdS) \cong 
D^b(\mdE)$ also work in the $\ZZ^n$-graded context. That is, 
the functors $\bL : D^b(\mmEn) \to D^b(\mmSn)$ and 
$\bR : D^b(\mmSn) \to D^b(\mmEn)$ are defined by the same way as 
the $\ZZ$-graded case, and they 
give an equivalence $D^b(\mmSn) \cong D^b(\mmEn)$, see 
\cite[Theorem~4.1]{Y1}. Note that the dualizing complex $\D$ of $S$ is 
$\ZZ^n$-graded, and $\DS(-) = \Hom^\bullet_S(-, \D)$ is also 
a duality functor on $D^b(\mmSn)$. Similarly, $\DE(-) = \Hom_E(-,E)$ 
is a duality functor on $D^b(\mmEn)$. As Proposition~\ref{dualities}, 
for $\N \in D^b(\mmEn)$, we have $\DS \circ \bL (\N) \cong 
\bL \circ \DE (\N) (-{\bf 2})[2n]$ in $D^b(\mmSn)$. 
Here we set ${\bf j} := (j,j,\ldots, j) \in \NN^n$ for $j \in \ZZ$. 

For $\ba = (a_1, \ldots, a_n)\in \ZZ^n$, 
set $\supp (\ba) := \{i  \mid a_i > 0\} \subset [n] := \{1, \ldots, n \}$ 
and $|\ba| = \sum_{i=1}^n a_i$. 
We say  $\ba \in \ZZ^n$ is {\it squarefree} if $a_i= 0,1$ for all $i \in [n]$. 
When $\ba \in \ZZ^n$ is squarefree, we sometimes identify $\ba$ with  
$\supp (\ba)$. 

\begin{dfn}[\cite{Y}]
We say a $\ZZ^n$-graded $S$-module $M$ is {\it squarefree}, if 
the following conditions are satisfied. 
\begin{itemize}
\item[(a)] $M$ is $\NN^n$-graded 
(i.e., $M_\ba =0$ if $\ba \not \in \NN^n$) and finitely generated. 
\item[(b)] The multiplication map 
$M_{\ba} \ni y \mapsto (\prod x_i^{b_i}) \cdot y \in M_{\ba +\bb}$ is 
bijective for all $\ba, \bb \in \NN^n$ with $\supp (\ba+\bb) = \supp (\ba)$. 
\end{itemize}
\end{dfn}

For a simplicial complex $\Delta \subset 2^{[n]}$, the Stanley-Reisner ideal 
$I_\Delta := ( \, \prod_{i \in F} x_i \mid F \not \in \Delta \, )$ 
and the Stanley-Reisner ring  $K[\Delta] := S/I_\Delta$ 
are squarefree modules. 
Note that if $M$ is squarefree then $M_\ba \cong M_F$ as $K$-vector
spaces for all $\ba \in \NN^n$ with $\supp(\ba)=F$. 
Let $\Sq$ be the full subcategory of $\mmSn$ consisting 
of squarefree modules.  In $\mmSn$, $\Sq$ is closed under 
kernels, cokernels and  extensions (\cite[Lemma~2.3]{Y}), and we have that 
$D^b(\Sq) \cong D^b_{\Sq}(\mmSn)$. If $\M \in D^b(\Sq)$, then 
$\DS(\M) \in D^b_{\Sq}(\mmSn)$ (see \cite{Y1}). 
So $\DS$ gives a duality functor on $D^b(\Sq)$.

\begin{dfn}[R\"omer~\cite{R0}]
A $\ZZ^n$-graded $E$-module $N = \bigoplus_{\ba \in \ZZ^n} N_\ba$ is 
{\it squarefree} if $N$ is finitely generated and 
$N = \bigoplus_{F \subset [n]} N_{-F}$. 
\end{dfn}

For example, any monomial ideal of $E$ is a squarefree module. 
Any monomial ideal of $E$ is of the form 
$J_\Delta = ( \, \prod_{i \in F} y_i \mid  F \not \in \Delta \, )$ 
for some simplicial complex $\Delta \subset 2^{[n]}$. We say 
$K\{\Delta\} := E/J_\Delta$ is the {\it exterior face ring} of $\Delta$. 

Let $\SqE$ be the full subcategory of $\mmEn$ consisting of squarefree 
$E$-modules. Then there exist functors $\cS: \SqE \to \Sq$ 
and $\cE: \Sq \to \SqE$ giving an equivalence $\Sq \cong \SqE$. Here 
$\cS(N)_F = N_{-F}$ for $N \in \SqE$, and the multiplication map 
$\cS(N)_F \ni z \mapsto x_i z \in \cS(N)_{F \cup \{ i \}}$ for $i \not \in F$ 
is given by  $\cS(N)_F =N_{-F} \ni z \mapsto (-1)^{\alpha(i, F)} y_i z 
\in N_{-(F \cup \{ i \})}= \cS(N)_{F \cup \{ i \}}$, where 
$\alpha(i,F) =  \# \{ \, j \in F \mid j < i \, \}$. 
For example, $\cS(K\{\Delta\}) = K[\Delta]$.  
See \cite{R0} for further information. Of course, $\cS$ and $\cE$ 
can be extended to the functors between $D^b(\Sq)$ and $D^b(\SqE)$. 

If $N \in \SqE$, then $N^* = \Hom_E(N,E)$ is squarefree again. 
So $(-)^*$ gives the duality functor $\DE$ on $D^b(\SqE)$. 
For example, $K\{\Delta\}^* = J_{\Delta^\vee}$, where $\Delta^\vee 
= \{ F \subset [n] \mid [n] \setminus F \,  
\not \in \Delta \}$ is the Alexander dual complex of $\Delta$. We have the 
{\it Alexander duality functor} $\bA := \cS \circ \DE \circ \cE$ on 
$\Sq$ (or $D^b(\Sq)$). Of course, $\bA(K[\Delta]) = I_{\Delta^\vee}$. 
In general, we have $\bA(H^i(\M))_F = (H^{-i}(\M)_{[n] \setminus F})^\vee$.

An associated prime ideal of $M \in \mmSn$ is of the form 
$P_F = ( x_i \mid i \not \in F )$ for some $F \subset [n]$. 
Let $M \in \Sq$ be a squarefree module.  
A monomial prime ideal $P_F$ is a minimal prime of $M$ 
if and only if $F$ is a maximal element of 
the set $\{ G \subset [n] \mid M_G \ne 0 \}$. 
The following is a squarefree version of Definition~\ref{dist pair for S}. 

\begin{dfn}\label{dist pair for SqS}
We say $(F, i) \in 2^{[n]} \times \ZZ$ is a {\it distinguished pair} for 
a complex $\M \in D^b(\Sq)$, if $P_F$ is a minimal prime of 
$H^i(\M)$ and $H^j(\M)_G =0$ for all $j$ with $j < i$ and 
$G \supset F$ with $|G| < |F| + i-j$. 
\end{dfn}

\begin{thm}\label{main2}
Let $\M \in D^b(\Sq)$. A pair $(F, i)$ is distinguished for $\M$ 
if and only if $(F, -|F|-i)$ is distinguished for $\DS(\M)$. 
If this is the case, $\dim_K H^i(\M)_F = 
\dim_K H^{-|F|-i}(\DS(\M))_F$. 
\end{thm}

\begin{proof}
Like the proof of Theorem~\ref{main}, 
we consider the spectral sequence $E_2^{p,q} = \Ext^p_S(H^{-q}(\M), \D) 
\Rightarrow \Ext^{p+q}_S(\M, \D)$. 
Then $E_r^{p,q}$ is squarefree for all $p,q$ and 
$r \geq 2$. When we consider a distinguished pair $(F, i)$, we set 
$$\dim_F M := 
\begin{cases}
-\infty & \text{if $M_G = 0$ for all $G \supset F$} \\
\max \{ \, |G| \mid G \supset F, \, M_G \ne 0 \, \} & 
\text{otherwise} 
\end{cases}
$$
for $M \in \Sq$. Set  $d_i(\M):= \dim_F H^i(\M)$ 
and $e_2^{p,q} := \dim_F \Ext^p_S(H^{-q}(\M), \D)$ 
for $\M \in D^b(\Sq)$. We also remark that 
$\dim_K M_F = l_{S_{P_F}}(M \otimes_S S_{P_F})$ for $M \in \Sq$. 
The equation \eqref{basic} 
holds in this context, and the proof of Theorem~\ref{main} works verbatim. 
\end{proof} 

If $\N \in D^b(\SqE)$, then 
it is easy to see that $\bL(\N)(-\b1) \in D^b(\Sq)$. 
So $\cL(-) := \bL(-)(-\b1)$ gives a functor from $D^b(\SqE)$ to $D^b(\Sq)$.
Moreover, we have $\cL \cong \bA \circ \DS \circ \cS$ by 
\cite[Proposition~4.3]{Y1}. 

\begin{dfn}
Let $\N \in D^b(\SqE)$. We say $(F,i)$ is a {\it distinguished pair} for $\N$
if it is a distinguished pair for $\cL(\N) \in D^b(\Sq)$ in the sense of 
Definition~\ref{dist pair for SqS}. 
\end{dfn}

The next result can be proved by the same way as Corollary~\ref{dist dual} 
using Theorem~\ref{main2}. 

\begin{prop}\label{dist dual for SqE}
Let $\N \in D^b(\SqE)$. A pair $(F,i)$ is distinguished for $\N$ 
if and only if $(F, 2n-|F|-i)$ is distinguished for $\DE(\N)$. 
If this is the case, we have $\dim_K H^i(\cL(\N))_F = 
\dim_K H^{2n-|F|-i}(\cL \circ \DE(\N))_F$. 
\end{prop}

If $\M \in D^b(\mmSn)$, then 
$\Tor^S_i(K,\M) := H^{-i}(K \otimes_E \P)$ is a $\ZZ^n$-graded module, 
where $\P$ is a graded free resolution of $\M$. 
Set $\beta_{i,\ba}(\M) := \dim_K \Tor^S_i(K,\M)_\ba$ for $\ba \in \ZZ^n$. 
We say $\beta_{i,\ba}(\M)$ is the $(i,\ba)$th Betti number of $\M$. 
If $\M \in D^b(\Sq)$ and $\beta_{i,\ba}(\M) \ne 0$, then $\ba$ is squarefree 
(see \cite{Y1}). 
 
\begin{dfn}[c.f. {\cite{BCP}}]
A Betti number $\beta_{i,F}(\M) \ne 0$ is {\it extremal} 
if $\beta_{j,G}(\M) = 0$ for all $(j,G) \ne (i,F)$ with 
$j \geq i$, $G \supset F$, and $|G|-j > |F|-i$.
\end{dfn}

Some of known results and backgrounds of extremal Betti numbers 
are found in the introduction of the present paper.

\begin{prop}[c.f. {\cite{AH}}]\label{dist <-> extremal}
Let $\M \in D^b(\Sq)$ and $\N := \cE(\M) \in D^b(\SqE)$. 
A pair $(F,i)$ is distinguished for $\DE(\N)$ if and only if 
$\beta_{i+|F|-n,F}(\M)$ is an extremal Betti number. 
If this is the case, then  
$\beta_{i+|F|-n,F}(\M) = \dim_K H^i(\cL \circ \DE(\N))_F$. 
\end{prop}

\begin{proof}
For $j \in \ZZ$ and $G \subset [n]$, we have the following. 
\begin{eqnarray*}
\beta_{j, G}(\M) 
&=& \dim_K [H^{|G|-j-n}(\DS \circ \bA (\M))]_{[n] \setminus G} 
\quad \text{(by \cite[Corollary~3.6]{Y1})}\\ 
&=& \dim_K [H^{n+j-|G|}(\bA \circ \DS \circ \bA (\M))]_G \\
&=& \dim_K [H^{n+j-|G|}(\cL \circ \cE \circ \bA (\M))]_G \\
&=& \dim_K [H^{n+j-|G|}(\cL \circ \DE (\N))]_G.        
\end{eqnarray*}
The assertion easily follows from this equality. 
\end{proof}

\begin{cor}\label{BCP for complexes}
Let $\M \in D^b(\Sq)$. A Betti number $\beta_{i,F}(\M)$ is extremal 
if and only if so is $\beta_{|F|-i, F}(\bA(\M))$. 
If this is the case,  
$\beta_{i,F}(\M) = \beta_{|F|-i, F}(\bA(\M))$. 
\end{cor}

\begin{proof}
If $\beta_{i,F}(\M)$ is extremal, then $(F, n+i-|F|)$ is a distinguished 
pair for $\DE \circ \cE (\M)$ by Proposition~\ref{dist <-> extremal}. 
By Proposition~\ref{dist dual for SqE}, 
$(F, n-i)$ is a distinguished pair for 
$\cE(\M) \cong \DE \circ \cE \circ \bA(\M)$. 
So $\beta_{|F|-i,F}(\bA(\M))$ is extremal. 
The converse implication can be proved by the same way. 
The last equality follows from Proposition~\ref{dist dual for SqE}. 
\end{proof}

This corollary generalizes results of Bayer-Charalambous-Popescu~\cite{BCP}, 
R\"omer~\cite{R0} and Miller~\cite{Mil}. 
Roughly speaking, the above proof is a ``complex version" 
of \cite{R0}. But, his argument itself does not work in the 
$\ZZ^n$-graded context, since he use a generic base change of $V = E_{-1}$.

For $\M \in D^b(\Sq)$. Set 
$\pd(\M) = \max\{ \, i \mid \text{$\beta_{i,F}(\M) \ne 0$ for some $F$} \, \}$ 
and $\reg(\M) = \max\{ \, |F|-i \mid \beta_{i,F}(\M) \ne 0 \, \}$. 
Since Betti numbers $\beta_{i,F}(\M)$ which give $\pd(\M)$ or $\reg(\M)$ 
are extremal, the next result follows from Corollary~\ref{BCP for complexes}. 

\begin{cor}[c.f.{\cite{Mil, R0}}]
If $\M \in D^b(\Sq)$, then $\pd(\M) = \reg(\bA(\M))$. 
\end{cor}


\begin{thebibliography}{99}
\bibitem{AAH} A. Aramova, L. Avramov and J. Herzog, 
Resolutions of monomial ideals and cohomology over exterior algebras. 
Trans. Amer. Math. Soc. {\bf 352} (2000), 579-594.

\bibitem{AH}  A. Aramova and J. Herzog, 
Almost regular sequences and Betti numbers, 
Amer. J. Math. {\bf 122} (2000), 689-719.

\bibitem{BCP} D. Bayer, H. Charalambous and S. Popescu, 
Extremal Betti numbers and applications to monomial ideals, J. Algebra 
{\bf 221} (1999), 497--512. 

\bibitem{BTN} E. Babson, I. Novik and  R. Thomas, 
Symmetric iterated Betti numbers, preprint (math.CO/0206063).

\bibitem{BH} W. Bruns and J. Herzog,  Cohen-Macaulay rings, revised edition,  
Cambridge University Press, 1998.

\bibitem{EFS} D. Eisenbud, G. Fl\o ystad  and F.-O. Schreyer, 
Sheaf cohomology and free resolutions over exterior algebra, 
Trans. Amer. Math. Soc. {\bf 355} (2003), 4397-4426.  

\bibitem{EPY} D. Eisenbud, S. Popescu and S. Yuzvinsky, 
Hyperplane arrangement cohomology and monomials in the exterior algebra, 
Trans. Amer. Math. Soc. {\bf 355} (2003), 4365-4383.   

\bibitem{GM} S.I. Gelfand and Y.I. Manin, 
Methods of homological algebra. 
Springer-Verlag, 1996. 

\bibitem{Kalai} 
G. Kalai, Algebraic shifting,  in 
``Computational commutative algebra and combinatorics" (T. Hibi, ed.), 
pp. 121--163, Adv. Stud. Pure Math. 33, Math. Soc. Japan, 2002.    

\bibitem{Mil} E. Miller, 
The Alexander duality functors and local duality with monomial support, 
J. Algebra. {\bf 231} (2000), 180--234. 

\bibitem{R0} T. R\"omer,  Generalized Alexander duality and applications, 
Osaka J. Math. {\bf 38} (2001), 469--485.  

\bibitem{Y} K. Yanagawa, 
Alexander duality for Stanley-Reisner rings and squarefree 
$\NN^n$-graded modules, J. Algebra {\bf 225} (2000), 630--645.  

\bibitem{Y1} K. Yanagawa, 
Derived category of squarefree modules and local cohomology with monomial 
ideal support, to appear in J. Math. Soc. Japan. 


\end{thebibliography}
\end{document}